\documentclass[12pt]{article}

\usepackage{amsmath, amsfonts, amssymb}
\usepackage{mathrsfs}
\usepackage{theorem}
\usepackage{bm}

\newtheorem{thm}{Theorem}[section]

\newtheorem{lem}[thm]{Lemma}

\newtheorem{defn}[thm]{Definition}

{\theorembodyfont{\rmfamily} \newtheorem{rem}[thm]{Remark}}
{\theorembodyfont{\rmfamily} }
\newcommand{\norm}[1]{\lVert#1\rVert_{r}^{2}}

\def\R{\mathbb R}

\newcommand{\fin}{\hspace*{\fill}\rule{0.3em}{1ex}}

\numberwithin{equation}{section}

\allowdisplaybreaks

\begin{document}

\title{Stability in Distribution  of Neutral Stochastic Functional Differential Equations  with Infinite Delay
}

\author{Hussein K. Asker ${}^{a,b}$ \footnote{ Email:husseink.askar@uokufa.edu.iq }  \\[0.2cm]
{\small a: Department of Mathematics, Faculty of Computer Science and Mathematics, }\\
{\small Kufa University, Al-Najaf, Iraq }\\
{\small b:Department of Mathematics, Swansea University, Singleton Park, SA2 8PP, UK }\\
 }
\maketitle
\begin{abstract}
In this paper, we investigate stability in distribution of neutral stochastic functional differential equations with infinite delay (NSFDEwID) at the state space\\
 \begin{equation*}
 C_{r}=\{{\varphi\in C((-\infty,0];R^{d}):\|\varphi\|_{r}=\sup_{-\infty<\theta\leq0}e^{r\theta}\lvert\varphi(\theta)\rvert} < \infty, \quad r > 0 \}.
 \end{equation*}
 We drive a sufficient strong monotone condition for the existence and uniqueness of the global solutions of NSFDEwID in the state space $ C_{r} $. We also address  the stability of  the solution map $ x_{t} $  and illustrate the theory with an example.
\end{abstract}
\section{Introduction}
 Recently, the asymptotic properties of neutral stochastic functional differential equations (NSFDEs) have been addressed by many authors. Kolmanovskii and Myshkis  \cite{KM} have  studied the theory of existence and uniqueness,  and  the stability of  the solutions of  these equations. Mao  \cite{M1, M2} has investigated  the moment exponential stability and the almost sure exponential stability.  Bao and Hou \cite{BH}, Under a non-Lipschitz condition and a weakened linear growth condition, the existence and uniqueness of mild solutions to neutral partial functional SDEs have been investigated. Tan et al. \cite{TJS} by the weak convergence approach, they have reviewed stability in distribution for NSFDEs. Huang and Deng  \cite{HD} have studied Razumikhin-type theorems on general $ p $-th moment asymptotic stability. Kolmanovskii et al. \cite{KMM} have dealt with NSDDE with Markovian switching. Liu and Xia \cite{LX} have established some results which are more effective and relatively easy to verify to obtain the required stability, Hu and Wang \cite{hw} have concerned with the NSFDEs with Markovian switching and derived some sufficient conditions for stability in distribution, to name a few. On the other hand,  infinite NSFDEs have also been  received a great attention, many remarkable results got on the existence and uniqueness and the asymptotic behaviour of NSFDEwID.  Ren and Xia \cite{RX}, for example, have studied existence and uniqueness of solutions with infinite delay at phase space $ BC((- \infty; 0]; \mathbb{R}^{d}) $ which denotes the family of bounded continuous  $ \mathbb{R}^{d} $-value functions  with norm $ \lVert \varphi\lVert = \sup_{ - \infty < \theta \le 0} \rvert \varphi (\theta)\rvert $ under non-Lipschitz condition. Prato and Zabczyk \cite{dz} devoted their work to asymptotic properties of solution. Boufoussi and Hajji \cite{BH2, BH1}, by using successive approximation, they have proved existence and uniqueness result for a class of NFSDEs in Hilbert spaces. Bao \cite{Bao} has considered the existence and uniqueness of solutions in $ L^{p}(\Omega, C_{h}) $ space. Caraballo et al. \cite{CRT} have analysed the almost sure exponential stability and ultimate boundedness of the solutions. For Wei and Cai \cite{WC} by choosing $ C_{g} $ as  phase space, the existence and uniqueness of the solutions to neutral stochastic functional differential equations with infinite delay have obtained under non-Lipschitz condition, weakened liner growth condition and contractive condition. Chen and Banas\cite{CB},  have devoted their work to obtain some sufficient conditions for the exponential stability as well as almost surely exponential stability for mild solution of neutral stochastic partial differential equations with delays by establishing an integral-inequality. Zhou, Y. \cite{Z} has studied the stability property of stochastic differential equations in Hilbert spaces. Bao and Cao \cite{BC} have examined the existence and uniqueness of mild solutions. Bao et al. \cite{BYY, BYYW} have investigated the ergodicity for several kinds of FSDEs with variable delays. and Wu et al. \cite{FYM} have been studied the ergodicity of underlying processes and establishes existence of the invariant measure for SFDEs with infinite delay.\\
The state space for functional differential equations (FDEs) with infinite delay plays a key role in solving a specific problem. In general, the state space will be Banach space of functions or equivalence class of functions see \cite{KS}. Mohammed \cite{M} has examined solution maps of SFDEs with finite delay on appropriate phase spaces and proved that the solution maps have Markov property. Based on the Markov property of solution maps of SFDEs with finite delay, Bao et al.\cite{BYY, BYYW} have examined the ergodicity, while Wu et al. \cite{FYM} investigates existence and uniqueness of solutions, Markov properties, and ergodicity of SFDEs with infinite delay by using phase state $ C_{r} $.\\
Motivated by the discussion above, in this paper, to guarantee the wellposedness of solutions of NSFDEwID, boundedness of solutions  from different initial data and further asymptotic properties including the mean-square boundedness and convergence of the solutions from different initial data we have imposed appropriate strong monotone condition. We have also studied, the stability in distribution of the solution data. Finally, we have introduced an example to describe the theory that we addressed.
\section{Perliminary}
Throughout this paper, unless otherwise specified, we use the following notation. Let  $ R^{d} $ denote the usual $d$-dimensional Euclidean space and $ \lvert \cdotp \lvert $   the Euclidean norm. If $A$ is a vector or a matrix, its transpose is denoted by $ A^{T} $; and $  \lvert A \rvert = \sqrt{\mbox{trace} (A^{T}A)} $ denotes its trace norm. Denote by $ x^{T}y $ the inner product of  $x$ and $y $ in $\R^{d} $. Let  $ C((-\infty,0];R^{d}) $ denote the family of all  continuous functions from $ (-\infty,0] $  to $ R^{d}.$ We choose the state space with the fading memory  defined as follows: for given positive number $ r $,
\begin{equation}\label{2.1}
C _{r} = \Big \{ \varphi \in C((-\infty,0];R^{d}): \|\varphi \| _{r} = \sup_{ - \infty < \theta \leq 0} e^{r\theta} \arrowvert \varphi (\theta) \arrowvert < \infty  \Big \}.
\end{equation}
 $ C_{r} $ is a Banach space with norm $ \|\cdot \| _{r}$, which is introduced in  \cite{Kuang and Smith}, contains the Banach space of bounded and continuous functions and for any $0<r_{1}\leq r_{2} < \infty$, $C_{r_{1}} \subset C_{r_{2}} $. Let $(\Omega,\mathcal{F}, \mathbb{P}) $  be a complete probability space with a filtration $ {\{{\mathcal{F}_{t}\}_{t\in [0,+\infty)}}} $ satisfying the usual conditions (i.e. it is right continuous and $ {{\mathcal{F}_{0}}} $  contains all P-null sets). Let $  1_{B} $ denote the indicator function of a set B. $  M^{2}([ a, b];R^{d})$ is a family of process $ \{ x(t) \}_{a \leq t \leq b}  $ in $ \mathcal{L}^{2}([a, b]; \mathbb{R}^{d}) $ such that $ \mathbb{E} \int_{a}^{b} \rvert x(t)\rvert^{2} dt < \infty $. Consider a $ d $-dimensional neutral stochastic functional differential equations with infinite delay
\begin{equation}\label{2.2}
d\{x(t) - D(x_{t})\} = b(x_{t})dt + \sigma (x_{t})dw(t), \quad \text{on} \quad t \geq 0,
\end{equation}
with the initial data
\begin{equation}\label{2.3}
x_{0} = \xi = \{\xi(\theta):-\infty<\theta\leq0 \} \in C_{r},
\end{equation}
where
\begin{equation*}
x_{t} = \{ x(t+\theta): - \infty < \theta \leq 0\}
\end{equation*}
 and $ b, D:  C_{r} \to\mathbb R^{d} $;      $ \sigma : C_{r} \to\mathbb R^{d\times m} $ are Borel measurable, $ w(t) $  is an $ m $-dimensional Brownian motion. It should be pointed out that $ x(t) \in  R^{d}  $ is a point, while $ x_{t}  \in C_{r}$
 is a continuous function on the interval $ (-\infty, t] $ taking values in $ R^{d}. $ Now, we  give the definition of the solutions for the equation \eqref{2.2}.
\begin{defn}\label{d1}
    \cite{FYM,Mo} Set $ {\mathcal{F}_{t}} = {\mathcal{F}_{0}} $ for $ -\infty< t \leq 0 $ and let x(t),   $ -\infty< t \leq \tau _{e}$ be continuous $ {R}^{d} $-value and $ {\mathcal{F}_{t}} $-adapted process.  $ x(t) $ is called a $   \textbf {local  strong  solution}   $ of \eqref{2.2} with initial data  $ \xi \in C_{r} $ if $ x(t)= \xi(t) \quad on - \infty < t \leq 0 $ and for all $ t \geq 0 $,
    \begin{equation*}
    x (t \wedge \tau _{k}) =   \xi (0) + D ( x _{t} \wedge \tau _{k} ) - D ( \xi) + \int_{0}^{t \wedge \tau _{k}} b (x_{s}) ds + \int_{0}^{t \wedge \tau _{k}} \sigma ( x_{s}) dw(s) \quad a.s.,
    \end{equation*}
    for each $ \ell \geq 1 $, where $ \{ \tau _{\ell} \} _{\ell \geq 1} $  is a non-decreasing sequence of stopping times such that  $ \tau _{\ell} \longrightarrow \tau _{e} $ almost surely as $ \ell \longrightarrow \infty $. If moreover, $ \lim \sup_{t \longrightarrow \tau _{e} } \arrowvert x(t) \arrowvert  = \infty   $  is satisfied a.s. when $ \tau _{e} < \infty  $ a.s., x(t) $ ( - \infty < t < \tau _{e} ) $ is called a \textbf{maximal local strong solution}  and $  \tau _{e} $ is called the explosion time. It is called a global solution when $ \tau _{e} = \infty $. A maximal local strong solution x(t), $ - \infty < t < \tau _{e} $ is said to be unique if for any other maximal local strong solution $ \overline{x} (t) $, $ - \infty < t < \overline {\tau } _{e} $, we have $ \tau _{e} = \overline {\tau } _{e} $ and $ x(t) = \overline{x}(t) $ for $ - \infty < t < \tau _{e} $ almost surely.
\end{defn}
\section{Existence and uniqueness of solutions}
In this section we study the existence and uniqueness of the global solutions of the NSFDEwID \eqref{2.2}, mean-square boundedness and convergence of the solutions from different initial data. Let $ \mathcal{P}(C_{r} )  $ denotes the family of all probability measures on $ (C_{r},\mathcal{B}(C_{r} )) $. Denote $ \mathcal{C}_{b}(C_{r}) $ the set of all bounded continuous functional. For any $ F \in \mathcal{C}_{b}(C_{r}) $, $ F:C_{r} \rightarrow \mathbb{R} $  and $ \pi(\cdot) \in  \mathcal{P}(C_{r} ) $, let $ \pi(F ):=\int_{C_{r}} F(\phi) \pi(d\phi)$. $ M_{0} $ stands for the set of probability measures  on $ (- \infty, 0] $, namely, for any $ \mu \in M_{0} $, $ \int_{- \infty}^{0} \mu (d\theta) = 1  $. For any $ r > 0 $, let us further define $ {M}_{r}  $  as follows:
\begin{equation}
M _{r} : = \Big\{ \mu \in M_{0}; \mu ^{(r)}:=\int_{- \infty}^{0} e^{-r\theta} \mu (d\theta) < \infty \Big\}.
\end{equation}
Obviously, there exist many such probability measures, see \cite{FYM, huwuhuang} for details.
\begin{rem}\label{remark 1}
\cite{FYM} Fix $ r_{_0 } > 0 $ and $  \mu \in M _{r_{0}} $. For any $ r \in (0, r_{0}) $, $ \mu^{(r)} $ is continuously non-decreasing on $ r $ and satisfies $ \mu^{(r_{0})} > \mu^{(r)} > \mu^{(0)} = 1$, and $ M_{r_{0}} \subset M_{r} \subset M_{0} $.
\end{rem}
In order to examine the stability in distribution of the equation \eqref{2.2} we assume the following conditions.
\begin{description}
\item[(H1)]  For $ \mu_{1} \in M _{2r}$ and $ \varphi \in C_{r} $ there exist $ k \in (0, 1)  $ with $  \mu^{(2r)}_{1} < 1 $  such that
\begin{equation}\label{4.6}
\rvert D (\varphi)-D(\phi) \rvert ^{2} \leq k \int_{- \infty}^{0} \rvert \varphi(\theta)-\phi(\theta) \rvert ^{2} \mu_{1} (d\theta) \quad \text{and} \quad D(0)=0.
\end{equation}
\item[(H2)]  Let $ b $ be a continuous functions. Assume there exist constants $ \lambda_{1}, \lambda_{2}, \lambda_{3}, \lambda_{4}>0 $, and probability measures $ \mu_{2}, \mu_{3} \in M_{2r} $ such that for any $ \varphi , \phi \in C_{r} $
\begin{equation}\label{4.7}
\begin{split}
\Big[ \varphi (0) - \phi (0) - \big( D(\varphi) - D (\phi) \big) \Big]^{T} \Big[ b ( \varphi) - b ( \phi ) \Big] &\leq  - \lambda_{1} \rvert \varphi (0) - \phi (0) \rvert ^{2} \\
& + \lambda_{2} \int_{-\infty}^{0} \rvert \varphi (\theta) - \phi (\theta) \rvert ^{2} \mu _{2} (d \theta),
\end{split}
\end{equation}
and for any function $ \sigma $
\begin{equation}\label{4.8}
\big \rvert \sigma (\varphi) - \sigma (\phi) \big\rvert ^{2} \leq \lambda_{3} \big\rvert \varphi (0) - \phi (0) \big\rvert ^{2} + \lambda_{4} \int_{- \infty}^{0} \big \rvert \varphi (\theta) - \phi (\theta)   \big \rvert ^{2} \mu _{3} (d\theta).
\end{equation}
\end{description}
\begin{rem}
	For simplicity, we suppose that $ \mu_{1} = \mu_{2} = \mu_{3} = \mu $.
\end{rem}
\begin{rem}
	It is easy to get from \eqref{4.6}, \eqref{4.7} and \eqref{4.8}, that
	\begin{description}
		\item[i)]  \begin{equation}\label{i}
		\rvert D (\varphi) - D (\phi)\rvert ^{2} \leq k \mu^{(2r)}\norm{\varphi-\phi}.
		\end{equation}
		\item[ii)] \begin{equation}\label{ii}
		\begin{split}
		\Big[ \varphi (0) - \phi (0) - \big( D(\varphi) - D (\phi) \big) \Big]^{T} \Big[ b ( \varphi) - b ( \phi ) \Big] &\leq  (- \lambda_{1}+ \lambda_{2} \mu^{(2r)})\norm{\varphi-\phi}.
		\end{split}
		\end{equation}
		\item[iii)]		\begin{equation}\label{iii}
		\big \rvert \sigma (\varphi) - \sigma (\phi) \big\rvert ^{2} \leq (\lambda_{3}+ \lambda_{4} \mu^{(2r)})\norm{\varphi-\phi}.
		\end{equation}
	\end{description}
\end{rem}
So, the assumptions \noindent{\bf (H1)} and \noindent{\bf (H2)} guarantees that the equation \eqref{2.2} admits a unique local strong solution $ \{x(t;\xi)\}_{t>- \infty} $ with initial data $ \xi \in C_{r} $. So, for global solution and further asymptotic properties including the mean-square boundedness and convergence of the solutions from different initial data, we impose that conditions.
\begin{lem}\label{lemma 2}
	\cite{Mo}   Let $ p > 1 $, $ \varepsilon > 0 $ and $ a, b \in R $. Then $ \rvert a + b \rvert ^{p} \leq \big [ 1 + \varepsilon ^{\frac{1}{p-1}} \big]^{p-1}  \Big( \rvert a \rvert ^{p} + \dfrac{\rvert b \rvert ^{p}}{\varepsilon} \Big) $.
\end{lem}
\begin{rem}
When $ \phi \equiv 0, D(\phi) \equiv 0 $, $ \varepsilon_{1},\varepsilon_{2}>0 $, by the Lemma \ref{lemma 2} and the assumption \noindent{\bf (H2)}, we get the strong monotone condition:
\begin{equation} \label{monotone}
\begin{split}
2\Big[ \varphi (0)  -  D(\varphi)  \Big]^{T} \Big[ b ( \varphi)  \Big]  +  \big \rvert \sigma (\varphi)  \big\rvert ^{2}\leq  - \alpha_{1} \rvert \varphi (0)  \rvert ^{2}  + \alpha_{2} \int_{-\infty}^{0} \rvert \varphi (\theta)  \rvert ^{2} \mu  (d \theta) + N ,
\end{split}
\end{equation}
where, $ \alpha_{1} = 2 \lambda_{1} - 2 \varepsilon_{1} -\dfrac{\lambda_{3}}{1-\varepsilon_{2}}  $,
$ \alpha_{2} =  2 \lambda_{2} +  2k\varepsilon_{1} +\dfrac{\lambda_{4}}{1-\varepsilon_{2}} $, $ N =\dfrac{1}{\varepsilon}_{1} \rvert b(0) \rvert ^{2} + \dfrac{1}{\varepsilon}_{2}  \rvert \sigma (0)\rvert^{2}  $.
\end{rem}
\begin{lem}\label{l1} For any $ \xi \in C_{r} $, under the assumption $ \bf(H1) $ we have,
\begin{equation*}
\rvert \xi(0) -  D(\xi) ) \rvert^{2}\leq M  \rVert \xi \rVert^{2}_{r},
\end{equation*}
where $ M=(1+k) ( 1  + \mu^{(2r)})  $
\end{lem}
\noindent{\bf Proof:} By the Lemma \ref{lemma 2} and the assumption $ \bf(H1) $ we have,
    \begin{equation*}
    \begin{split}
    \rvert \xi(0) -  D(\xi) ) \rvert^{2}  &\leq (1+\varepsilon) \rvert \xi(0) \rvert^{2} + ( 1 +\dfrac{1}{\varepsilon} )  \rvert D(\xi) \rvert^{2}\\
    &  \leq (1+\varepsilon) \sup_{ - \infty < \theta \le 0} e^{2r\theta }\rvert \xi(\theta) \rvert^{2} + k ( 1 +\dfrac{1}{\varepsilon} )  \int_{- \infty}^{0}  e^{2r\theta} \rvert \xi(\theta) \rvert ^{2} e^{-2r\theta} \mu (d\theta)\\
    &  \leq (1+\varepsilon) \sup_{ - \infty < \theta \le 0} e^{2r\theta }\rvert \xi(\theta) \rvert^{2} + k ( 1 +\dfrac{1}{\varepsilon} )  \int_{- \infty}^{0} \sup_{ - \infty < \theta \le 0} e^{2r\theta} \rvert \xi(\theta) \rvert ^{2} e^{-2r\theta} \mu (d\theta)\\
    & \leq (1+\varepsilon) \rVert \xi \rVert^{2}_{r} + k ( 1 +\dfrac{1}{\varepsilon} ) \int_{- \infty}^{0} \rVert \xi \rVert^{2} e^{-2r\theta} \mu (d\theta) \\
    & = (1+\varepsilon) \rVert \xi \rVert^{2}_{r} + k ( 1 +\dfrac{1}{\varepsilon} )  \rVert \xi \rVert^{2}  \mu^{(2r)}.\\
    \end{split}
    \end{equation*}
Set $ \varepsilon=k $, so that 
 \begin{equation}\label{ln1}
\begin{split}
\rvert \xi(0) -  D(\xi) ) \rvert^{2}  &\leq (1+k) ( 1  + \mu^{(2r)})    \rVert \xi \rVert^{2}_{r}=M\rVert \xi \rVert^{2}_{r}.\\
\end{split}
\end{equation}
 \begin{lem}\label{l2} Under $ \bf(H1) $ let $ \xi,\eta \in C_{r} $, then
\begin{equation}\label{l6}
\rvert \xi(0) - \eta(0) - ( D(\xi) - D ( \eta )) \rvert^{2}  \leq M \rVert \xi - \eta \rVert_{r}^{2},
\end{equation}
where, $ M=(1+k) ( 1  + \mu^{(2r)}) $.\\
The proof is similar to that of Lemma \ref{l1}, we omit it here.
 \end{lem}
\begin{lem}\label{Lf1} Let $\xi\in
	C_r((-\infty,0];\R^d)$ and $ \bf(H1) $ holds, then for $ t>0 $ there exist positive constants $k_1$
	and $k_2$ such that
	\begin{equation}\label{lc1}
	\begin{split}
	\sup_{0<s\leq t}|\xi(s)|^2&\leq
	k_1e^{-2rs}\|\xi\|^2_r+k_2\sup_{0<s\leq t}|\xi(s)-D(\xi_s)|^2,
	\end{split}
	\end{equation}
	where $k_1=\frac{k
		\mu^{(2r)}}{1-k}$ and $k_2=\frac{1}{(1-k)^{2}}$.
\end{lem}
\noindent{\bf Proof:}
By the lemma \ref{lemma 2} and condition $ \bf(H1) $, for any
$\varepsilon>0$,
\begin{equation*}
\begin{split}
&|\xi(t)|^2=|D(\xi_t)+\xi(t)-D(\xi_t)|^2\\
&\leq(1+\varepsilon)\Big(\frac{|D(\xi_t)|^2}{\varepsilon}+|\xi(t)-D(\xi_t)|^2\Big)\\
&\leq(1+\varepsilon)\Big(\frac{k^2}{\varepsilon}\int_{-\infty}^0|\xi(t+\theta)|^2\mu
d(\theta) +|\xi(t)-D(\xi_t)|^2\Big).
\end{split}\end{equation*}
Taking $\varepsilon=\frac{k}{1-k}$ and using the definition of norm in
the phase space $C_r((-\infty,0];\R^d)$, we obtain
\begin{equation*}
\begin{split}
\sup_{0<s\leq t}|\xi(s)|^2 &\leq k\int_{-\infty}^0\sup_{0<s\leq
	t}|\xi(s+\theta)|^2\mu (d\theta)+\frac{1}{(1-k)}\sup_{0<s\leq
	t}|\xi(s)-D(\xi_s)|^2\\
& =k\int_{-\infty}^{-s}\sup_{0<s\leq	t}|\xi(s+\theta)|^2\mu (d\theta) +k\int_{-s}^0\sup_{0<s\leq	t}|\xi(s+\theta)|^2\mu (d\theta) \\
&+\frac{1}{(1-k)}\sup_{0<s\leq	t}|\xi(s)-D(\xi_s)|^2\\
& \leq k\mu^{(2r)} e^{-2rs}\|\xi\|_r^2 +k\sup_{0<s\leq t}|\xi(s)|^2+\frac{1}{(1-k)}\sup_{0<s\leq t}|\xi(s)-D(\xi_s)|^2.
\end{split}\end{equation*}
Subsequently,
\begin{equation*}
\begin{split}
\sup_{0<s\leq t}|\xi(s)|^2 &\leq \frac{k\mu^{(2r)}}{1-k}
e^{-2rs}\|\xi\|_r^2 +\frac{1}{(1-k)^{2}}\sup_{0<s\leq
	t}|\xi(s)-D(\xi_s)|^2\\
& = k_1 e^{-2rs}\|\xi\|_r^2 +k_2\sup_{0<s\leq
	t}|\xi(s)-D(\xi_s)|^2,\\
\end{split}
\end{equation*}
where $k_1=\frac{k \mu^{(2r)}}{1-k}$ and
$k_2=\frac{1}{(1-k)^{2}}$. The proof is complete.\hfill $\Box$\\
\begin{lem}\label{Lf2} Let  $r>0$ and $\xi, \eta\in
	C_r((-\infty,0];\R^d)$. Let condition $ \bf(H1) $ hold. Then there
	exist positive constants $k_3, k_{4}$, such that
  \begin{equation*}
\sup_{0\leq s\leq t}|\xi(s)-\eta(s)|^2 \leq  k_3e^{-2rs}\|\xi - \eta\|^2_r+
k_4\sup_{0\leq s\leq t}|\xi(s)-\eta(s)-D(z_s)+D(\eta_s)|^2,
  \end{equation*}
  where $k_3=\frac{k\mu^{(2r)}}{(1-k)}$ and
  $k_4=\frac{1}{(1-k)^{2}}$.
\end{lem}
Since the proof is similar to that of Lemma \ref{Lf2}, we will not give details here.
\begin{rem}
In the next theorem, by the initial data $ \xi $, we mean the initial function or initial segment process.  So that $ \xi $ is a function not a fixed constant. To highlight the initial segment process, we denote by $ x(t; \xi) $ and $ x_{t}(\xi) $  the solution and the solution map of \eqref{2.2}, respectively.
\end{rem}
\begin{thm}\label{4.}
    Under assumptions  $ \bf(H1) $ and $ \bf(H2) $,
    \begin{description}
        \item[(i)] For any initial data $  \xi \in C_{r}  $ the NSFDEwID \eqref{2.2} has a global solution $ x(t) $ almost surely, which is continuous and $ \mathcal{F}_{t} $-adapted.
        \item[(ii)]  If $ \lambda_{1}, \lambda_{2}, \lambda_{3}$ and $ \lambda_{4} $ satisfy $ 2 \lambda_{1} > 73 \lambda_{3}  + 2 \lambda_{2} \mu ^{(2r)} +73 \lambda_{4} \mu ^{(2r)} $, then there exist constants $ C_{1}, C_{2} > 0 $ and $ \lambda \in \big( 0,  \dfrac{1}{M}\big[2 \lambda_{1} - 73 \lambda_{3} - 2 \lambda_{2} \mu ^{(2r)} -73 \lambda_{4}\mu ^{(2r)}\big]\wedge 2r \big) $ such that for any initial data $ \xi \in C_{r} $, 
        \begin{equation}
        \mathbb{E} \big(  \rvert x(t;\xi)  \rvert ^{2} \big) \leq C_{1} + C_{2}\mathbb{E}  \rVert \xi \rVert^{2}_{r} e^{-\lambda t},
        \end{equation}
where $ C_{1}= \dfrac{2k_{2}}{\lambda} \Big( \dfrac{73 \rvert \sigma (0) \rvert^{2}}{\varepsilon_{2}}   +  \dfrac{\rvert b(0) \rvert^{2}}{\varepsilon_{1}} \Big) $,\\$ C_{2}=\Bigg\{k_{1}+2k_{2} \Bigg[M + \dfrac{\mu^{(2r)}}{2r - \lambda} \Big(  (1 +k) \lambda+ 2 \lambda_{2}+ 2k\varepsilon_{1} +  \dfrac{73\lambda_{4}}{1-\varepsilon_{2}} \Big)   \Bigg] \Bigg\} $ and $ \varepsilon_{1}, \varepsilon_{2} $ are both sufficiently small constants such that
\begin{equation*}
 \Bigg[2 \lambda_{1}- M\lambda - 2\varepsilon_{1}- \dfrac{73 \lambda_{3}}{1-\varepsilon_{2}} -\Big(  2 \lambda_{2}+2k\varepsilon_{1}+ \dfrac{73 \lambda_{4}}{1-\varepsilon_{2}}  \Big) \mu^{(2r)}     \Bigg] > 0,
\end{equation*}
namely, solution $ x(t;\xi)  $   is mean-square bounded.
        \item[(iii)] Under conditions in $ \bf(ii) $, for different initial dat $ \xi $ and $ \eta $, the corresponding solution $ x(t; \xi) $ and $ x(t;\eta) $ satisfy
         \begin{equation}
         \mathbb{E} \Big( \sup_{ 0 < s \le t}  \rvert  x(t;\xi) - x(t; \eta) \rvert^{2} \Big)  \leq C_{3} \mathbb{E}\norm{\xi -\eta}e^{-\lambda t},
         \end{equation}
        where  $ C_{3}=\Bigg\{k_{3}+2k_{4} \Bigg[ M + \dfrac{\mu^{(2r)}}{2r - \lambda} \Big(  (1 +k) \lambda+ 2 \lambda_{2} + 73 \lambda_{4} )       \Big)                \Bigg] \Bigg\}. $
    \end{description}
\end{thm}
\noindent{\bf Proof:} We divide the proof into three steps:\\
\noindent{\bf Step 1 ( Proof of (i).):}
The proof is similar to that of \cite[Theorem 3.2]{FYM}, we here
only highlight the difficulty from the neutral term. By the conditions  $ \bf(H1) $ and $ \bf(H2) $, the coefficients are local Lipschitz continuous, then there exists a uinque maximal local strong solution $ x(t) $ to \eqref{2.2} on $ t \in (- \infty, \tau_{\ell}) $. In order to prove
that $x(t)$ is a global solution,  we need only prove that $
\tau_{e} = \infty $ almost surely. Define $  \tau_{\ell} = \inf \big\{ t \geq 0 : \rvert x(t) \rvert \geq \ell
\big\}  $, then $ \tau_{\ell} $ is increasing as $ \ell  \rightarrow  \infty $, let $ \lim\limits_{\ell \rightarrow \infty} \tau_{\ell}=
 \tau_{\infty} \leq \tau_{e} $ almost surly. If we can
show $ \tau_{\infty} = \infty $ a.s., then  $ \tau_{e} = \infty $
a.s., which implies that $ x(t) $ is actually global. This is
equivalent to prove that for any $ T > 0 $, $ \mathbb{P} (
\tau_{\ell} \leq T) \rightarrow 0 $ as $ \ell \rightarrow \infty $.\\
Using $ \bf(H1) $, one can derive that
\begin{equation}\label{4.16}
\begin{split}
&\sup_{ 0 <  s \le t \wedge \tau_{\ell}}  \Big (\mathbb{E} \rvert  x(s) \rvert^{2} \Big)  \leq ( 1 + \varepsilon)  \sup_{ 0 <  s \le t \wedge \tau_{\ell}}  \Big (  \mathbb{E} \rvert  x(s)  - D   ( x_{s})\rvert ^{2} \Big)\\
& \quad + k ( 1+ \dfrac{1}{\varepsilon}) \sup_{ 0 <  s \le t \wedge \tau_{\ell}}  \Big (\mathbb{E}\int_{- \infty}^{0}  \rvert x(s+\theta) \rvert ^{2}  \mu (d\theta) \Big)\\
& \leq ( 1 + \varepsilon)  \sup_{ 0 <  s \le t \wedge \tau_{\ell}}  \Big (  \mathbb{E} \rvert  x(s)  - D   ( x_{s})\rvert ^{2} \Big)+ k ( 1+ \dfrac{1}{\varepsilon}) \mathbb{E} \Big[\int_{- \infty}^{-s}  \sup_{ 0 <  s \le t \wedge \tau_{\ell}}  \Big ( \rvert x(s+\theta) \rvert ^{2} \Big)  \mu (d\theta) \\
&\quad + \int_{- s}^{0}  \sup_{ 0 <  s \le t \wedge \tau_{\ell}}  \Big ( \rvert x(s+\theta) \rvert ^{2} \Big)  \mu (d\theta)\Big]\\
&= ( 1 + \varepsilon)  \sup_{ 0 <  s \le t \wedge \tau_{\ell}}  \Big (  \mathbb{E} \rvert  x(s)  - D   ( x_{s})\rvert ^{2} \Big)\\
& \quad+ k ( 1+ \dfrac{1}{\varepsilon}) \mathbb{E} \Big[\int_{- \infty}^{-s}  \sup_{ 0 <  s \le t \wedge \tau_{\ell}}  \Big (e^{2r(s+\theta)} \rvert x(s+\theta) \rvert ^{2} \Big) e^{-2r(s+\theta)} \mu (d\theta) \\
&\quad + \int_{- s}^{0}  \sup_{ 0 <  s \le t \wedge \tau_{\ell}}  \Big ( \rvert x(s+\theta) \rvert ^{2} \Big)  \mu (d\theta)\Big]\\
&\leq ( 1 + \varepsilon)  \sup_{ 0 <  s \le t \wedge \tau_{\ell}}  \Big (  \mathbb{E} \rvert  x(s)  - D   ( x_{s})\rvert ^{2} \Big)\\ 
&+ k ( 1+ \dfrac{1}{\varepsilon}) \mathbb{E} \Big[\sup_{ -\infty <  s+\theta \le t \wedge \tau_{\ell}-s}  \Big (e^{2r(s+\theta)} \rvert x(s+\theta) \rvert ^{2} \Big)\int_{- \infty}^{-s}   e^{-2r(s+\theta)} \mu (d\theta) \\
&\quad + \sup_{ -\theta <  s+\theta \le t \wedge \tau_{\ell}}  \Big ( \rvert x(s+\theta) \rvert ^{2} \Big)\int_{- s}^{0}    \mu (d\theta)\Big]\\
&\leq ( 1 + \varepsilon)  \sup_{ 0 <  s \le t \wedge \tau_{\ell}}  \Big (  \mathbb{E} \rvert  x(s)  - D   ( x_{s})\rvert ^{2} \Big)\\
&\quad+ k ( 1+ \dfrac{1}{\varepsilon}) \mathbb{E} \Big[\sup_{ -\infty <  s+\theta \le 0}  \Big (e^{2r(s+\theta)} \rvert x(s+\theta) \rvert ^{2} \Big)\int_{- \infty}^{-s}   e^{-2r(s+\theta)} \mu (d\theta)\\
& \quad + \sup_{ 0 <  s \le t \wedge \tau_{\ell}}  \Big ( \rvert x(s) \rvert ^{2} \Big)\int_{- \infty}^{0}    \mu (d\theta)\Big]\\
&= ( 1 + \varepsilon)  \sup_{ 0 <  s \le t \wedge \tau_{\ell}}  \Big (  \mathbb{E} \rvert  x(s)  - D   ( x_{s})\rvert ^{2} \Big)+ k e^{-2rs} ( 1+ \dfrac{1}{\varepsilon}) \mu^{(2r)}\mathbb{E} \lVert \xi \lVert _{r}^{2}  \\
& \quad+ k ( 1+ \dfrac{1}{\varepsilon}) \sup_{ 0  <  s \le t \wedge \tau_{\ell}} \Big( \mathbb{E}  \rvert x(s) \rvert ^{2}   \Big) . 
\end{split}
\end{equation}
Choosing $ \varepsilon > \dfrac{k}{1-k} $  implies $ \gamma=k ( 1+ \dfrac{1}{\varepsilon})<1 $, we arrive at  
\begin{equation}\label{4.17}
\begin{split}
\sup_{ 0 <  s \le t \wedge \tau_{\ell}}  \Big (\mathbb{E} \rvert & x(s) \rvert^{2} \Big)  \leq  \dfrac{1 + \varepsilon}{(1- \gamma )}  \sup_{ 0 <  s \le t \wedge \tau_{\ell}}  \Big (  \mathbb{E} \rvert  x(s)  - D   ( x_{s})\rvert ^{2} \Big) + \dfrac{\gamma e^{-2rs}  \mu^{(2r)}}{1- \gamma}\mathbb{E} \rVert \xi \rVert^{2}_{r}.
\end{split}
\end{equation}
Applying the It\^o  formula to $ |x(t) - D(x_{t})|^{2} $, by
\eqref{ln1} and the monotone condition \eqref{monotone} yields for
any $ t \in [0, T] $,
\begin{equation}\label{s1}
\begin{split}
\mathbb{E} \rvert & x(t \wedge \tau_{\ell}) - D   ( x_{t \wedge \tau_{\ell}})\rvert^{2} = \mathbb{E} \rvert x(0) - D ( \xi) \rvert^{2} + 2 \mathbb{E} \int_{0}^{t \wedge \tau_{\ell}} [ x(s) - D (x_{s})]^{T} b( x_{s }) ds\\
& \qquad + \mathbb{E} \int_{0}^{t \wedge \tau_{\ell}} \rvert \sigma ( x_{s}) \rvert ^{2} ds\\
 &  \leq M \mathbb{E} \rVert \xi \rVert^{2}_{r} -\alpha_{1}   \mathbb{E} \int_{0}^{t \wedge \tau_{\ell}} \rvert x (s)  \rvert ^{2} ds + \alpha_{2} \mathbb{E} \int_{0}^{t \wedge \tau_{\ell}} \int_{-\infty}^{0} \rvert x (s+\theta) \rvert ^{2} \mu (d \theta) ds  +  N .
\end{split}
\end{equation}
By the fact [(3.12), from \cite{FYM}], we have
\begin{equation}\label{u}
\begin{split}
&\int_{0}^{t \wedge \tau_{\ell}} \int_{- \infty}^{0}   \rvert x(s+\theta) \rvert ^{2} \mu (d\theta) ds \leq \dfrac{1}{2r} \lVert \xi \lVert _{r} ^{2} \mu ^{(2r)}  +  \int_{0}^{t }  \rvert x(s \wedge \tau_{\ell}) \rvert ^{2}  ds.
\end{split}
\end{equation}
 Substituting \eqref{u} into \eqref{s1}, one has
\begin{equation}\label{4.14}
\begin{split}
\mathbb{E} \rvert x(t \wedge \tau_{\ell}) - D ( x_{t \wedge \tau_{\ell}})\rvert^{2}  &  \leq   L_{1} + J_{1}  \int_{0}^{t } \mathbb{E} \rvert x(s \wedge \tau_{\ell}) \rvert ^{2}  ds,
\end{split}
\end{equation}
where, $ L_{1} =\Big[ \Big(M + \dfrac{\alpha_{2}\mu^{(2r)}}{2r}  \Big) \mathbb{E} \rVert \xi \rVert^{2}_{r} + N \Big] $ and $  J_{1} =  - \alpha_{1}+\alpha_{2} $.\\
Therefore, by substituting   \eqref{4.14} into \eqref{4.17} we obtain
\begin{equation}\label{4.18}
\begin{split}
\mathbb{E} \rvert x(t \wedge \tau_{\ell}) \rvert^{2}  &  \leq  \dfrac{(1 + \varepsilon)L_{1}}{(1- \gamma )}  +   \dfrac{J_{1}(1 + \varepsilon)}{(1- \gamma)}   \mathbb{E} \int_{0}^{t}    \rvert  x(s \wedge \tau_{\ell} ) \rvert ^{2} ds \\
& +\dfrac{\gamma e^{-2rs}  \mu^{(2r)}}{1- \gamma} \mathbb{E} \rVert \xi \rVert^{2}_{r}. 
\end{split}
\end{equation}
Hence, by applying  Gronwall's inequality yields
\begin{equation}\label{4.19}
\begin{split}
\mathbb{E} \rvert x(t \wedge \tau_{\ell}) \rvert^{2}  &  \leq   L_{2}      e^{J_{2}t} , \\
\end{split}
\end{equation}
where, $ L_{2} = \dfrac{(1 + \varepsilon)L_{1}}{(1- \gamma )}  +\dfrac{\gamma e^{-2rs}  \mu^{(2r)}}{1- \gamma} \mathbb{E} \rVert \xi \rVert^{2}_{r} $, $ J_{2} = \dfrac{J_{1}(1 + \varepsilon)}{(1- \gamma )} . $ 

According to the definition of $ \tau_{\ell}, $ we have
\begin{equation}
\begin{split}
\mathbb{E} \Big (   \rvert x( T \wedge \tau_{\ell}) \rvert^{2} \Big)&  =    \mathbb{E} \Big ( \rvert x( T \wedge \tau_{\ell}) \rvert^{2}  I_{\{ T \leq \tau_{\ell} \}}\Big) + \mathbb{E} \Big ( \rvert x( T \wedge \tau_{\ell}) \rvert^{2}  I_{\{ T \geq \tau_{\ell} \}}\Big)\\  
& \ge \ell^{2} \mathbb{P} \big( \tau_{\ell} \leq T \big).
\end{split}
\end{equation}
Consequently, by \eqref{4.19}
\begin{equation*}
\begin{split}
\mathbb{P} ( \tau_{\ell} \leq T) & \leq  \dfrac{L_{2}      e^{J_{2}t}}{\ell^{2}}.
\end{split}
\end{equation*}
This implies
\begin{equation*}
\limsup_{\ell \rightarrow \infty} \mathbb{P} ( \tau_{\ell} \leq T ) = 0,
\end{equation*}
which means that \eqref{2.2} has a unique global solution $ x(t) $ on $ [0, \infty) $ almost surely.\\

\noindent{\bf Step 2 ( Proof of (ii).):} For any $ \lambda \in \big( 0,  \dfrac{1}{M}\big[2 \lambda_{1} - 73 \lambda_{3} - 2 \lambda_{2} \mu ^{(2r)} -73 \lambda_{4}\mu ^{(2r)}\big]\wedge 2r \big) $ applying the It\^o  formula to $ e ^{\lambda t}|x(t) - D(x_{t})|^{2}, $   together with the  \eqref{ln1}, the Lemma \ref{lemma 2} with $ \varepsilon = k $, the monotone condition \eqref{monotone}  and the assumption  \noindent{\bf (H1)},     yields for any $ t \in [0, T] $ that,
\begin{equation}\label{4.42}
\begin{split}
&\mathbb{E}  \Big ( \sup_{ 0 < s\le t}  e ^{\lambda s} \rvert x(s)  - D ( x_{s})\rvert^{2} \Big)  = \mathbb{E} \rvert x(0) - D ( \xi) \rvert^{2} +  \mathbb{E} \Big ( \sup_{ 0 <  s \le t}\int_{0}^{s} e ^{\lambda u} \Big [ \lambda \rvert x(u) - D (x_{u}) \rvert^{2} \\
&  + 2[x(u) - D (x_{u})] ^{T} b( x_{u }) + \rvert \sigma (x_{u}) \rvert^{2} \Big] du \Big) + 2\mathbb{E} \Big ( \sup_{ 0 <  s \le t}\int_{0}^{s} e ^{\lambda u}  [x(u) - D (x_{u})] ^{T} \sigma (x_{u}) dw(u) \Big) \\
& \leq M \mathbb{E}  \rVert \xi \rVert^{2}_{r} + \lambda  (1+k) \mathbb{E} \Big ( \sup_{ 0 <  s \le t} \int_{0}^{s} e^{\lambda u}  \rvert x(u)\rvert ^{2} du \Big) +  \lambda(1 + \dfrac{1}{k})  \mathbb{E} \Big( \sup_{ 0 <  s \le t} \int_{0}^{s} e^{\lambda u}  \rvert D (x_{u}) \rvert^{2} du \Big)\\
&   +  \mathbb{E} \Big ( \sup_{ 0 <  s \le t} \int_{0}^{s} e ^{\lambda u} \Big[ - \alpha_{1} \big\rvert x (u)  \big\rvert ^{2} + \alpha _{2} \int_{- \infty}^{0} \big \rvert x (u+\theta) \big \rvert ^{2} \mu  (d\theta) +  N \Big] du \Big)\\
& + 2\mathbb{E} \Big ( \sup_{ 0 <  s \le t}\int_{0}^{s} e ^{\lambda u}  [x(u) - D (x_{u})] ^{T} \sigma (x_{u}) dw(u) \Big) \\
& \leq M \mathbb{E}  \rVert \xi \rVert^{2}_{r} +   ((1+k)\lambda-\alpha_{1}) \mathbb{E} \Big ( \sup_{ 0 <  s \le t} \int_{0}^{s} e^{\lambda u}  \rvert x(u)\rvert ^{2} du \Big)\\
&   + ((1+k)\lambda+\alpha_{2}) \mathbb{E} \Big ( \sup_{ 0 <  s \le t} \int_{0}^{s} e ^{\lambda u}   \int_{- \infty}^{0} \big \rvert x (u+\theta) \big \rvert ^{2} \mu  (d\theta)du\Big) + \dfrac{ N}{\lambda}(e^{\lambda t}-1)  \\
& + 2\mathbb{E} \Big ( \sup_{ 0 <  s \le t}\int_{0}^{s} e ^{\lambda u}  [x(u) - D (x_{u})] ^{T} \sigma (x_{u}) dw(u) \Big). \\
\end{split}
\end{equation}
By [(3.13), \cite{FYM}], we have
\begin{equation}\label{4.24}
\begin{split}
\int_{0}^{t} \int_{- \infty}^{0} e^{\lambda s} & \rvert x(s +\theta) \rvert ^{2} \mu (d\theta) ds  
\leq \dfrac{1}{2r - \lambda} \rVert \xi \rVert^{2} _{r}  \mu ^{(2r)} + \mu  ^{(2r)} \int_{0}^{t } e^{\lambda s} \rvert x(s) \rvert ^{2}  ds .
\end{split}
\end{equation}
Now, using the Burkholder-Davis-Gundy inequality, the Lemma \ref{lemma 2}, \noindent{\bf (H1)}  and the condition \eqref{4.8},  we obtain
\begin{equation}\label{w}
\begin{split}
2\mathbb{E} \Big ( & \sup_{ 0 <  s \le t}  \int_{0}^{s} e ^{\lambda u}  [x(u) - D (x_{u})] ^{T} \sigma (x_{u}) dw(u) \Big)  \\
& \leq 8\sqrt{2} \mathbb{E} \Big ( \int_{0}^{t} e ^{2 \lambda s} \rvert [x(s) - D (x_{s})] ^{T} \sigma (x_{s})\rvert^{2} ds \Big) ^{\frac{1}{2}}\\
& \leq 12 \mathbb{E} \Big ( \int_{0}^{t} e ^{2 \lambda s} \rvert x(s) - D (x_{s}) \rvert ^{2}  \rvert \sigma (x_{s})\rvert^{2} ds \Big)^{\frac{1}{2}}\\
& \leq 12 \mathbb{E} \Big[ \Big ( \sup_{ 0 <  s \le t}  e ^{ \lambda s} \rvert x(s) - D (x_{s}) \rvert ^{2} \Big)^{\frac{1}{2}} \Big( \int_{0}^{t} e ^{ \lambda s} \rvert \sigma (x_{s})\rvert^{2} ds \Big)^{\frac{1}{2}} \Big]\\
&\leq \dfrac{1}{2}  \mathbb{E} \Big ( \sup_{ 0 <  s \le t}  e ^{ \lambda s } \rvert x(s) - D (x_{s}) \rvert ^{2} \Big) + 72 \mathbb{E}  \int_{0}^{t} e ^{ \lambda s} \rvert \sigma (x_{s})\rvert^{2} ds\\
& \leq \dfrac{1}{2} \mathbb{E} \Big ( \sup_{ 0 <  s \le t}  e ^{ \lambda s } \rvert x(s) - D (x_{s}) \rvert ^{2} \Big)+ \dfrac{72 \lambda_{3}}{1-\varepsilon_{2}} \mathbb{E}  \int_{0}^{t} e ^{ \lambda s}  \big\rvert x (s)  \big\rvert ^{2} ds\\
&  +  \dfrac{72  \lambda_{4}}{1-\varepsilon_{2}} \mathbb{E}  \int_{0}^{t} \int_{-\infty}^{0}  e ^{ \lambda s} \big\rvert x ( s + \theta)  \big\rvert ^{2} \mu  (d \theta) ds  +  \dfrac{72 }{ \lambda \varepsilon_{2} } \rvert \sigma (0)\rvert^{2} ( e^{\lambda t} - 1)   .
\end{split}
\end{equation}
By taking in account, the value of $ N, \alpha_{1}$ and $\alpha_{2} $,  substitute \eqref{4.24} and \eqref{w}  into \eqref{4.42} with $ \lambda \leq 2r$, we arrive at
\begin{equation}\label{4.25}
	\begin{split}
		\mathbb{E}& \Big ( \sup_{ 0 < s\le t}  e ^{\lambda s}  \rvert x(s) - D ( x_{s})\rvert^{2} \Big)  \leq 2 \Bigg[ M + \dfrac{\mu^{(2r)}}{2r - \lambda} \Big(  (1 +k) \lambda+ 2 \lambda_{2}+ 2k\varepsilon_{1} +  \dfrac{73\lambda_{4}}{1-\varepsilon_{2}}        \Big)           \Bigg]  \mathbb{E}  \rVert \xi \rVert^{2}_{r}\\
		& +  \dfrac{2}{\lambda} \Big( \dfrac{73 \rvert \sigma (0) \rvert^{2}}{\varepsilon_{2}}   +  \dfrac{\rvert b(0) \rvert^{2}}{\varepsilon_{1}} \Big) e^{\lambda t}  + 2 \Bigg[(1+k)\lambda - 2 \lambda_{1}+2\varepsilon_{1}+ \dfrac{73 \lambda_{3}}{1-\varepsilon_{2}} \\
		 & +\Big( (1+k)\lambda + 2 \lambda_{2}+2k\varepsilon_{1}+ \dfrac{73 \lambda_{4}}{1-\varepsilon_{2}}  \Big) \mu^{(2r)}   \Bigg] \mathbb{E} \Big ( \sup_{ 0 <  s \le t} \int_{0}^{s} e ^{\lambda u} \big\rvert x (u)  \big\rvert ^{2} du \Big)  .\\
	\end{split}
\end{equation}
Consequently, by the Lemma \ref{Lf1}, we have
\begin{equation}\label{4.26}
\begin{split}
\mathbb{E}\Big ( \sup_{ 0 < s\le t}&  e ^{\lambda s}  \rvert x(s) \rvert^{2} \Big)  \leq\Bigg\{k_{1}+2k_{2} \Bigg[M + \dfrac{\mu^{(2r)}}{2r - \lambda} \Big(  (1 +k) \lambda+ 2 \lambda_{2}+ 2k\varepsilon_{1} +  \dfrac{73\lambda_{4}}{1-\varepsilon_{2}} \Big)                    \Bigg] \Bigg\} \mathbb{E}  \rVert \xi \rVert^{2}_{r}\\
& +  \dfrac{2k_{2}}{\lambda} \Big( \dfrac{73 \rvert \sigma (0) \rvert^{2}}{\varepsilon_{2}}   +  \dfrac{\rvert b(0) \rvert^{2}}{\varepsilon_{1}} \Big) e^{\lambda t}  - 2k_{2} \Bigg[ 2 \lambda_{1}-M\lambda  -2\varepsilon_{1}- \dfrac{73 \lambda_{3}}{1-\varepsilon_{2}} \\
& -\Big(  2 \lambda_{2}+2k\varepsilon_{1}+ \dfrac{73 \lambda_{4}}{1-\varepsilon_{2}}  \Big) \mu^{(2r)}     \Bigg] \mathbb{E} \Big ( \sup_{ 0 <  s \le t} \int_{0}^{s} e ^{\lambda u} \big\rvert x (u)  \big\rvert ^{2} du \Big)  .
\end{split}
\end{equation}
Note that $ 2 \lambda_{1} > 73 \lambda_{3}  + 2 \lambda_{2} \mu ^{(2r)} +73 \lambda_{4} \mu ^{(2r)} $,  $ k\in(0, 1) $ and
\begin{equation*}
  \lambda \in \big( 0,  \dfrac{1}{M}\big[2 \lambda_{1} - 73 \lambda_{3} - 2 \lambda_{2} \mu ^{(2r)} -73 \lambda_{4}\mu ^{(2r)}\big]\wedge 2r \big),
\end{equation*} 
we can choose $ \varepsilon_{1}, \varepsilon_{2} $ sufficiently small such that:\\
$ \Bigg[2 \lambda_{1}- M\lambda - 2\varepsilon_{1}- \dfrac{73 \lambda_{3}}{1-\varepsilon_{2}} -\Big(  2 \lambda_{2}+2k\varepsilon_{1}+ \dfrac{73 \lambda_{4}}{1-\varepsilon_{2}}  \Big) \mu^{(2r)}     \Bigg] > 0 $, thus
\begin{equation}\label{4.27}
\begin{split}
\mathbb{E}& \Big ( \sup_{ 0 < s\le t}    \rvert x(s) \rvert^{2} \Big) \leq C_{1}+C_{2}\mathbb{E}  \rVert \xi \rVert^{2}_{r}e^{-\lambda t},\\
\end{split}
\end{equation}
where, 
$ C_{1}= \dfrac{2k_{2}}{\lambda} \Big( \dfrac{73 \rvert \sigma (0) \rvert^{2}}{\varepsilon_{2}}   +  \dfrac{\rvert b(0) \rvert^{2}}{\varepsilon_{1}} \Big) $ and \\$ C_{2}=\Bigg\{k_{1}+2k_{2} \Bigg[M + \dfrac{\mu^{(2r)}}{2r - \lambda} \Big(  (1 +k) \lambda+ 2 \lambda_{2}+ 2k\varepsilon_{1} +  \dfrac{73\lambda_{4}}{1-\varepsilon_{2}} \Big)                    \Bigg] \Bigg\} $.\\
This means that the solution $   x(t;\xi)  $ is mean-square bounded.\hfill $\Box$\\

\noindent{\bf Step 3 ( Proof of (iii).):} Let $  x(t; \xi) $ and $ x(t; \eta) $ be two different solutions with two different initial date  $ \xi, \eta $ to the equation \eqref{2.2}, then:
\begin{equation}
\begin{split}
d(x(t; \xi) - x(t; \eta) - D(x_{t }(\xi)) + D(x_{t }(\eta)) )& = \{b(x_{t }(\xi)) - b(x_{t }(\eta))\}dt \\
& + \{\sigma (x_{t }(\xi)) - \sigma (x_{t }(\eta))\}dw(t).
\end{split}
\end{equation}
For simplicity, let $ z(t) = x(t; \xi) - x(t; \eta) $ , $ \bar D(t) = D(x_{t }(\xi)) - D(x_{t }(\eta)) $ , $ \bar{b}(t) = b(x_{t }(\xi)) - b(x_{t }(\eta)) $ and $ \bar{\sigma}(t) = \sigma (x_{t }(\xi)) - \sigma (x_{t }(\eta)) $
with the initial data $ \xi - \eta $. For $ \lambda > 0  $ defined in (ii), applying the  It\^{o}   formula to $ e ^{\lambda t}\mathbb{E}\rvert z(t) -  \bar D(u)\rvert^{2} $ with \eqref{l6} gives:
\begin{equation}
\begin{split}
\mathbb{E} \Big( &\sup_{ 0 < s \le t} e^{ \lambda s} \rvert z ( s ) -  \bar D(s)\rvert^{2} \Big)  \leq M      \mathbb{E}  \rVert \xi - \eta \rVert_{r}^{2}\\
& \quad + \mathbb{E} \Big( \sup_{ 0 < s \le t} \int_{0}^{s} e^{\lambda u} \big[\lambda   \rvert z(u) - \bar D(u) \rvert ^{2} +2[z(u) -\bar D(u)  ] ^{T}  \bar{b}(z_{u}) + \rvert \bar{\sigma}(z_{u}) \rvert^{2}\big] du\\
& \quad + 2 \mathbb{E} \Big( \sup_{ 0 < s\le t} \int_{0}^{s} e^{\lambda u}[z(u) - \bar D(t) ] ^{T}  \bar{\sigma}(z_{u}) dw(u\Big).\\
\end{split}
\end{equation}
By the Lemma \ref{lemma 2} and the assumption $ \bf(H1) $, with $ \varepsilon=k $, we have
\begin{equation}
\rvert z(u) - \bar D(u) \rvert ^{2}\leq(1+k)   \rvert z(u)\lvert^{2}+(1+k)\int_{- \infty}^{0} \rvert z(u+\theta)\rvert^{2} \mu (d\theta),
\end{equation}
and by assumption $ \bf(H2) $,
\begin{equation}
2[z(u) -\bar D(u)  ] ^{T}  \bar{b}(z_{u})\leq - 2\lambda_{1}\lvert z(u)  \rvert^{2}+ 2 \lambda_{2}\int_{- \infty}^{0} \rvert z(u+\theta)\rvert^{2} \mu (d\theta),
\end{equation}
\begin{equation}
\begin{split}
&\lvert\bar{\sigma}(z_{u}) \rvert^{2} \leq \lambda_{3} \rvert z(u) \rvert^{2} +\lambda_{4} \int_{- \infty}^{0} \rvert z(u+\theta)\rvert^{2} \mu (d\theta),\\
\end{split}
\end{equation}
similar to \eqref{w}, with the condition \eqref{4.8}  we have
\begin{equation}
\begin{split}
&2 \mathbb{E} \Big( \sup_{ 0 < s\le t} \int_{0}^{s} e^{\lambda u}[z(u) - \bar D(t) ] ^{T}  \bar{\sigma}(z_{u}) dw(u\Big)  \leq\dfrac{1}{2}\mathbb{E} \Big( \sup_{ 0 < s \le t} e^{ \lambda s} \rvert z ( s ) -  \bar D(s)\rvert^{2} \Big)\\
& + 72\lambda_{3} \mathbb{E}\Big( \sup_{ 0 < s\le t}  \int_{0}^{s} e ^{ \lambda u}  \rvert z(u) \rvert^{2}du\Big) + 72\lambda_{4} \mathbb{E}\Big( \sup_{ 0 < s\le t}  \int_{0}^{s} \int_{- \infty}^{0}  e ^{ \lambda u}\rvert z(u+\theta)\rvert^{2} \mu (d\theta)du\Big).
\end{split}
\end{equation}
Consequently, by considering the fact \eqref{4.24} and the Lemma \ref{Lf2}, we get that
\begin{equation*}\label{4.30}
\begin{split}
\mathbb{E} &\Big( \sup_{ 0 < s \le t} e^{ \lambda s}  \rvert z ( s ) \rvert^{2} \Big)  \leq\Bigg\{k_{3}+2k_{4} \Bigg[ M + \dfrac{\mu^{(2r)}}{2r - \lambda} \Big(  (1 +k) \lambda+ 2 \lambda_{2} + 73 \lambda_{4} )       \Big)                \Bigg] \Bigg\}\mathbb{E}  \rVert \xi - \eta \rVert_{r}^{2} \\
& - 2k_{4} \Bigg[ 2 \lambda_{1}-73 \lambda_{3} -  M\lambda-2 \lambda_{2} \mu^{(2r)}  - 73 \lambda_{4}  \mu^{(2r)}   \Bigg]\mathbb{E}  \int_{0}^{t} e^{\lambda s}   \rvert z(s)\lvert^{2}ds .
\end{split}
\end{equation*}
Note that $ 2 \lambda_{1} > 73 \lambda_{3} +2 \lambda_{2} \mu^{(2r)}  + 73 \lambda_{4}  \mu^{(2r)}  $ and  \\ $ \lambda \in \big( 0,  \dfrac{1}{M}\big[2 \lambda_{1}-73 \lambda_{3} -2 \lambda_{2} \mu^{(2r)}  - 73 \lambda_{4}  \mu^{(2r)}\big] \wedge 2r \big) ,$  $ k\in(0, 1) $,  yields  that:\\
$ \Big[ 2 \lambda_{1}-73 \lambda_{3} -  M\lambda-2 \lambda_{2} \mu^{(2r)}  - 73 \lambda_{4}  \mu^{(2r)} \Big] > 0 $, implies to
\begin{equation}
\begin{split}
\mathbb{E}& \Big( \sup_{ 0 < s \le t}  \rvert  x(t;\xi) - x(t; \eta) \rvert^{2} \Big)  \leq C_{3} \mathbb{E}\norm{\xi -\eta}e^{-\lambda t}, 
\end{split}
\end{equation}
where, $ C_{3}=\Bigg\{k_{3}+2k_{4} \Bigg[ M + \dfrac{\mu^{(2r)}}{2r - \lambda} \Big(  (1 +k) \lambda+ 2 \lambda_{2} + 73 \lambda_{4} )       \Big)                \Bigg] \Bigg\} $.\\
The proof is therefore complete. \hfill $\Box$\\
\section{Stability in Distribution }
In this section, we shall study the stability in distribution for the  segment process  $ \{x_{ t}\}_{t \geq 0} $. We need the following Lemma.

\begin{lem}
		If the process $ \{x (t)\}_{t \geq 0} $ is the unique solution of the equation \eqref{2.2}, then the segment process $ \{ x_{t}\}_{t \geq 0 } $ is a strong homogeneous Markov process on $ C_{r} $:

	$ \mathbb{P}(x_{t} \in A) \mid \mathcal{F}_{s} ) = \mathbb{P}(x_{t} \in A) \mid x_{s} ) $ $ P $-a.s., for all $ 0 \leq s \leq t $ and Borel sets $  A \in \mathcal{B}(C_{r}) $.
\end{lem}

\noindent{\bf Proof:} We divide the proof into two steps ( the proof is similer to ( Theorem 4.2, \cite{FYM} )).

\textbf{Step 1: Strong Markov Property}: For all $ 0 \leq s \leq t < \infty $ and  a finite stopping time $ \tau > 0 $ be,   we consider the equation 
\begin{equation*}
\begin{split}
x (t + \tau) &= x(\tau) + \big( D ( x_{t +\tau}) - D (x_{\tau}) \big) + \int_{\tau}^{t } b(x_{s+\tau})ds + \int_{\tau}^{t} \sigma(x_{s+\tau})dw(s)\\
& =x(\tau) + \big( D ( x_{t +\tau}) - D (x_{\tau}) \big) + \int_{\tau}^{t+\tau } b(x_{s})ds + \int_{\tau}^{t+\tau} \sigma(x_{s})dw(s),
\end{split}
\end{equation*} 
By the definition of $ x_{t} $, we have $ x_{t}(t_{0}, x_{t_{0}})  = x ( t+\theta, t_{0} ; x(t_{0}+\theta)) $ with $ x_{t_{0}} (t_{0}, x_{t_{0}} ) = x_{t_{0}} $, for any $ t_{0}\leq t $. Note that $ w(t) $ is a strong Markov process with independent increment. It follows that $ \mathcal{F}_{\tau} $ is independent of $ \mathcal{G}_{\tau}=\sigma\{w(\tau+s) - w(\tau)\} $ for any $ s>0 $. Also, note that $ x_{t}(\tau, \xi) $ depends completely on the increments $ w(\tau+s) - w(\tau) $  and so is $ \mathcal{G}_{\tau}$-measurable when $ x_{\tau}=\xi $ is given. Hence $ x_{t}(\tau, \xi) $ is independent of $ \mathcal{F}_{\tau} $ for any $ t>0 $. For any $ A \subset C_{r} $, we therefore have
\begin{equation*}
\begin{split}
\mathbb{E} ( 1_{\{x_{t+\tau} (0, \xi)\}\in A} \mid \mathcal{F}_{\tau}) &= \mathbb{E} ( 1_{\{x_{t} (\tau, x_{\tau})\}\in A} \mid \mathcal{F}_{\tau}) \\
& = \mathbb{E} ( 1_{\{x_{t} (\tau, \xi)\}\in A}\rvert\xi=x_{\tau} )\\
&= \mathbb{P}(x_{\tau}=\xi, x_{t}(\tau, \xi)\in A\rvert\xi=x_{t})\\
&  = \mathbb{E} ( 1_{\{x_{t} (\tau, \xi)\}\in A} \mid x_{\tau})
\end{split}
\end{equation*}
Toward that standard technique, it takes after that for at whatever bounded Borel measurable function $ \varphi: C_{r} \rightarrow \mathbb{R} $,
\begin{equation*}
\mathbb{E} ( \varphi ( x_{t+\tau} ) \mid \mathcal{F}_{\tau} ] = \mathbb{E} ( \varphi ( x_{t+\tau} ) \mid x_{\tau} ),
\end{equation*}
which ends the proof and yields that $ x_{t} $ is a strong Markov process.\hfill $\Box$\\

\textbf{Step 2: Homogeneity:} According to the definition of transition probability:
\begin{equation*}
P( \xi , u; A, t+ u) = \mathbb{P} ( x _{t+u} ( u, \xi) \in A),
\end{equation*}
where $  x _{t+u} ( u, \xi)  $ is determined by the solution $ x(t) $:
\begin{equation*}
x( t+u) = \xi(0) + D ( x_{t +u}) - D(\xi) + \int_{u}^{t+u} b(x_{s}) ds + \int_{u}^{t+u} \sigma(x_{s}) dw(s).
\end{equation*}
This equation is equivalent to 
\begin{equation}\label{M-1}
x( t+u) = \xi(0) + D ( x_{t +u}) - D(\xi) + \int_{u}^{t} b(x_{s+u}) ds + \int_{u}^{t} \sigma(x_{s+u}) d\widetilde{w}(s),
\end{equation}
where $ \widetilde{w}(s) = w(s+u) - w(u) $ is clearly an $ d $-dimensional Brownian motion. So , we have
\begin{equation}\label{M-2}
x( t) = \xi(0) + D ( x_{t }) - D(\xi) + \int_{0}^{t} b(x_{s}) ds + \int_{0}^{t} \sigma(x_{s}) dw(s),
\end{equation}
with $  x_{0 } = \xi $. Comparing equations \eqref{M-1} with \eqref{M-2} and noting that $ x_{t+u} $ and $ x_{t} $ completely depends on $ x(t+u) $ and $ x(t) $ and their history, we see by the week uniqueness that $ \{ x_{t+u}\}_{t \geq 0}  $ are identical in probability law. Consequently,
\begin{equation*}
\mathbb{P}( x _{t+u} ( u, \xi) \in A ) = \mathbb{P} ( x_{t} ( 0, \xi)  \in A),
\end{equation*}
namely,
\begin{equation*}
P( \xi , u; A, t+ u) = P ( \xi, 0; A, t).
\end{equation*}
This complete the proof.\hfill $\Box$\\
\subsection{Stability in Distribution of  $ x_{t} $}
In order to address the stability of $ x_t $ in this subsection, let we highlight some notations, see \cite{TJS} for details. For  $ t\geq 0 $, denoted by $ P(\xi;t, \cdot) $ the transition probability of $ \{x_t\}_{t \geq 0}$ and for any $ P_{1}, P_{2}\in \mathcal{P}(C_{r}) $, define the metric $ d_{\mathbb{L}}  $ by 
\begin{equation*}
	d_{\mathbb{L}} (P_{1}, P_{2})=\sup_{ g \in \mathbb{L}}\Bigg| \int_{ C_{r} }g(\xi)(P_{1}(d\xi)- \int_{ C_{r} }g(\eta)(P_{2}(d\eta)\Bigg|,
\end{equation*}
where $ \mathbb{L}= \{g: C_{r}\rightarrow \mathbb{R}: \rvert g(\xi)-g(\eta) \rvert \leq \lVert \xi-\eta \lVert_{\infty}\quad \text{and}\quad \rvert g(\cdot)\rvert\leq1\quad \text{for}\quad \xi,\eta \in C_{r} \} $.
\begin{defn}
 \cite{TJS}. The process $ x_{t}(\xi) $ is said to be stable in distribution if there exists a probability measure $ \pi \in \mathcal{P}(C_{r}) $ such that $ P(\xi;t, \cdot) $ converges weakly to $ \pi $ as $ t \rightarrow \infty $ for any $ \xi \in C_{r} $, that is,
\begin{equation*}
\lim\limits_{t\rightarrow \infty}d_{\mathbb{L}}(P(\xi;t, \cdot),\pi(\cdot))=0, \quad \text{for all}\quad \xi \in  C_{r}.
\end{equation*}
In this case, \eqref{2.2} is said to be stable in distribution.
\end{defn}
We start by investigating some sufficient conditions on the stability in distribution for the segment process $ x_{t} $ on $ t\geq 0 $.
\begin{lem}\label{p1}
Under assumptions $ \bf(H1) $  and $ \bf(H2) $, if $ \lambda_{1}, \lambda_{2}, \lambda_{3}$ and $ \lambda_{4} $ satisfy $ 2 \lambda_{1} > 73 \lambda_{3}  + 2 \lambda_{2} \mu ^{(2r)} +73 \lambda_{4} \mu ^{(2r)} $, then there exist constants $ C_{4}, C_{5} > 0 $ and $ \lambda \in \big( 0,  \dfrac{1}{M}\big[2 \lambda_{1} - 73 \lambda_{3} - 2 \lambda_{2} \mu ^{(2r)} -73 \lambda_{4}\mu ^{(2r)}\big]\wedge 2r \big) $ such that for any initial data $ \xi \in C_{r} $, 
\begin{equation}\label{4.36}
\mathbb{E} \rVert x_{t} \rVert _{r}^{2} \leq C_{4} + C_{5} e ^{- \lambda t},
\end{equation}
where
\begin{equation*}
\begin{split}
& C_{4} = \dfrac{2k_{2}}{\lambda} \Big( \dfrac{73 \rvert \sigma (0) \rvert^{2}}{\varepsilon_{2}}   +  \dfrac{\rvert b(0) \rvert^{2}}{\varepsilon_{1}} \Big),\\
& C_{5} =\Bigg\{1+k_{1}+2k_{2} \Bigg[M + \dfrac{\mu^{(2r)}}{2r - \lambda} \Big(  (1 +k) \lambda+ 2 \lambda_{2}+ 2k\varepsilon_{1} +  \dfrac{73\lambda_{4}}{1-\varepsilon_{2}} \Big)  \Bigg] \Bigg\} \mathbb{E}  \rVert \xi \rVert^{2}_{r}.\\
\end{split}
\end{equation*}
\end{lem}
\noindent{\bf Proof:} Noting that $ \lambda < 2r $, therefore, correspond to the definition of the norm $ \rVert \cdot \rVert _{r} $, it is easy to see that:
\begin{equation}\label{5.8}
\begin{split}
\mathbb{E} \rVert x_{t} \rVert _{r}^{2} 
& =  e ^{- \lambda t} \mathbb{E} \rVert  \xi \rVert^{2}_{r} +  \mathbb{E} \Big( \sup_{ 0 < s \le t} \rvert x ( s ) \rvert^{2} \Big).
\end{split}
\end{equation}
From \eqref{4.27}, we have
\begin{equation}\label{5.4}
\begin{split}
\mathbb{E}& \Big ( \sup_{ 0 < s\le t} \rvert x(s) \rvert^{2} \Big)  \leq \dfrac{2k_{2}}{\lambda} \Big( \dfrac{73 \rvert \sigma (0) \rvert^{2}}{\varepsilon_{2}}   +  \dfrac{\rvert b(0) \rvert^{2}}{\varepsilon_{1}} \Big) \\
& + \Bigg\{k_{1}+2k_{2} \Bigg[M + \dfrac{\mu^{(2r)}}{2r - \lambda} \Big(  (1 +k) \lambda+ 2 \lambda_{2}+ 2k\varepsilon_{1} +  \dfrac{73\lambda_{4}}{1-\varepsilon_{2}} \Big)  \Bigg] \Bigg\}e^{-\lambda t} \mathbb{E}  \rVert \xi\rVert^{2}_{r}.\\
\end{split}
\end{equation}
By substituting \eqref{5.4} into \eqref{5.8}, we get
\begin{equation}\label{key}
\begin{split}
\mathbb{E} \rVert x_{t} \rVert _{r}^{2} & \leq C_{4} + C_{5} e ^{- \lambda t} ,\\
\end{split}
\end{equation}
 where
\begin{equation*}
\begin{split}
& C_{4} = \dfrac{2k_{2}}{\lambda} \Big( \dfrac{73 \rvert \sigma (0) \rvert^{2}}{\varepsilon_{2}}   +  \dfrac{\rvert b(0) \rvert^{2}}{\varepsilon_{1}} \Big),\\
& C_{5} = \Bigg\{1+k_{1}+2k_{2} \Bigg[M + \dfrac{\mu^{(2r)}}{2r - \lambda} \Big(  (1 +k) \lambda+ 2 \lambda_{2}+ 2k\varepsilon_{1} +  \dfrac{73\lambda_{4}}{1-\varepsilon_{2}} \Big)  \Bigg] \Bigg\} \mathbb{E}  \rVert \xi \rVert^{2}_{r}               .\\
\end{split}
\end{equation*}
This means that, the solution map $ x_{t} $ is mean-square bounded and giving the desired  assertion \eqref{4.36}.\hfill $\Box$
\begin{lem}\label{p2}
Under assumptions $ \bf(H1) $ and $ \bf(H2) $, if $ \lambda_{1}, \lambda_{2}, \lambda_{3}$ and $ \lambda_{4} $ satisfy $ 2 \lambda_{1} > 73 \lambda_{3}  + 2 \lambda_{2} \mu ^{(2r)} +73 \lambda_{4} \mu ^{(2r)} $, then there exist constant $ C_{6} > 0 $ and $ \lambda \in \big( 0,  \dfrac{1}{M}\big[2 \lambda_{1} - 73 \lambda_{3} - 2 \lambda_{2} \mu ^{(2r)} -73 \lambda_{4}\mu ^{(2r)}\big]\wedge 2r \big) $ such that
for the different initial data $ \xi $ and $ \eta  \in C_{r}$, the corresponding solution maps $  x_{t}(\xi) $ and $ x_{t}(\eta) $ satisfy,
\begin{equation}\label{5.2}
\mathbb{E} \rVert x_{t}(\xi) -  x_{t}(\eta) \rVert^{2} _{r} \leq C_{6} \mathbb{E} \rVert \xi -  \eta \rVert^{2} _{r} e^{-\lambda t},\\
\end{equation}
where
\begin{equation*}
C_{6} = \Bigg\{1+k_{3}+2k_{4} \Bigg[ M + \dfrac{\mu^{(2r)}}{2r - \lambda} \Big(  (1 +k) \lambda+ 2 \lambda_{2} + 73 \lambda_{4} )       \Big)                \Bigg] \Bigg\}.
\end{equation*}
Namely, the solution map $ x_{t} $ is mean-square bounded and solution maps from different initial data are convergent.
\end{lem}
Since the proof is comparable to that of Lemma \ref{p1}, we'll not explicate here.  

\begin{thm}
	Under assumptions $ \bf(H1) $ and $ \bf(H2) $, \eqref{2.2} is stable in distribution. 
\end{thm}
\noindent{\bf Prove:} Since \eqref{key} and \eqref{5.2} hold, the proof is standard, we omit here. 
\section{Example }
 In this section we give an example in order to verify the theory that we studied. Let we consider the one-dimensional neutral stochastic differential equation with infinity delay:
 \begin{equation}\label{ex}
 d[x(t) - D(x_{t})] = - c x(t) dt + \Big( cos(x(t))  + \int_{- \infty}^{0} x_{t} (\theta) \mu(d\theta) \Big) dw(t), \qquad t \geq 0,
 \end{equation}
 with initial value $ x(t) = \xi $ when $ t \in ( - \infty, 0] $. Where $ c > 0 $, $ w(t) $ is a scalar Brownian motion, while $ D : C_{r} \rightarrow R $ is defined by:
  \begin{equation*}
  D(\phi) = \dfrac{1}{2} \int_{ - \infty }^{0} \phi(\theta) \mu(d \theta).
  \end{equation*}
So that by  H\"{o}lder inequality we have
\begin{equation*}
\lvert D(\varphi) - D(\phi)\lvert ^{2} \leq \dfrac{1}{4} \int_{ - \infty }^{0}\lvert \varphi(\theta) - \phi(\theta) \lvert ^{2} \mu(d \theta).
\end{equation*}
Define $ b(\phi) = - c \phi, \sigma(\phi) = cos(\phi)  + \int_{- \infty}^{0} \phi (\theta) \mu(d\theta) $. It is easy to show that:
\begin{equation}\label{ex 1}
\begin{split}
\Big[ \varphi (0) - \phi (0) - \big( D(\varphi) - & D (\phi) \big) \Big]^{T} \Big[ b ( \varphi) - b ( \phi ) \Big] =   \Big[ \varphi (0) - \phi (0) \\
& \qquad + \dfrac{1}{2}\int_{ - \infty }^{0} \big(\varphi (\theta) - \phi (\theta )\big) \mu(d \theta) \Big]^{T} \Big[ - c\varphi  + c \phi\Big]\\
& \leq -c \rvert \varphi (0) - \phi (0) \rvert ^{2}  + \dfrac{c}{4} \int_{-\infty}^{0} \rvert \varphi (\theta) - \phi (\theta) \rvert ^{2} \mu  (d \theta),
\end{split}
\end{equation}
and
\begin{equation} \label{ex 2}
\begin{split}
\big \rvert \sigma (\varphi) - \sigma (\phi) \big\rvert ^{2} & =  \big\rvert cos(\varphi)  - cos(\phi)   +  \int_{- \infty}^{0} ( \varphi (\theta) - \phi (\theta)   ) \mu  (d\theta) \big\lvert^{2} \\
& \leq( 1+\varepsilon) \big\rvert cos(\varphi)  - cos(\phi) \big\rvert^{2} + \dfrac{1+\varepsilon}{\varepsilon} \big\rvert \int_{- \infty}^{0} ( \varphi (\theta) - \phi (\theta)   ) \mu  (d\theta) \big\lvert^{2}\\
& \leq (1+\varepsilon) \big\rvert \varphi(0)  - \phi(0) \big\rvert^{2} + \dfrac{1+\varepsilon}{\varepsilon}  \int_{- \infty}^{0} \big\rvert \varphi (\theta) - \phi (\theta)  \big\lvert^{2}  \mu  (d\theta). \\
\end{split}
\end{equation}
Thus, from \eqref{ex 1} and \eqref{ex 2} we get that: $ \lambda_{1}=c, \lambda_{2}=\dfrac{c}{4},\lambda_{3}=1+\varepsilon$ and $\lambda_{4}=\dfrac{1+\varepsilon}{\varepsilon} $, hence for any $ \varepsilon > 0 $ and $ \mu^{(2r)} \in M_{2r} $, the solution map $ x_{t} $ of the equation \eqref{ex} satisfies all the properties that we studied whenever $ c > \dfrac{4+ 146 (1+\varepsilon) + 146 (1+\varepsilon^{-1}) \mu^{(2r)}}{4 - \mu^{(2r)} } $.

\section*{Acknowledgement}
I would like to thank Professor Chenggui Yuan for his expert advice, as well as Dr. Jianhai Bao for his encouragement and advice throughout this project. This research was supported by Kufa University and the Iraqi Ministry of Higher Education and Scientific Research.

\end{document}